\documentclass{kluwer}    
\usepackage{amsfonts}
\usepackage{pstricks}
\usepackage{graphicx}

\newtheorem{theorem}{Theorem}

\newtheorem{definition}[theorem]{Definition}

\newtheorem{lemma}[theorem]{Lemma}
\newtheorem{proposition}[theorem]{Proposition}

\newenvironment{proof}[1][Proof]{\textbf{#1.} }{\ \rule{0.5em}{0.5em}}
\begin{document}
\begin{article}
\begin{opening}
\title{The Order Completion Method for Systems of Nonlinear PDEs:  Pseudo- topological Perspectives}
\author{Jan Harm \surname{van der Walt}}
\runningauthor{J H van der Walt} \runningtitle{The Order
Completion Method For Systems Of Nonlinear PDEs}
\institute{Department of Mathematics and Applied Mathematics\\
University of Pretoria}
\date{}

\begin{abstract}
By setting up appropriate uniform convergence structures, we are
able to reformulate the Order Completion Method of
Oberguggenberger and Rosinger in a setting that more closely
resembles the usual topological constructions for solving PDEs. As
an application, we obtain existence and uniqueness results for the
solutions of arbitrary continuous, nonlinear PDEs.
\end{abstract}


\classification{Mathematics Subject Clasifications (2000)} {34A34,
54A20, 06B30, 46E05}

\end{opening}

\section{Introduction and Preliminary Remarks}

\subsection{General, Type Independent Theories for PDEs}

It is widely believed that there is no general, type independent
theory concerning existence and regularity of solutions to
arbitrary nonlinear PDEs.  Indeed, the book `Lecture notes on
PDEs' by I. V. Arnold \cite{Arnold} starts with the following
remark:
\begin{eqnarray}
&&\mbox{``In contrast to ordinary differential equations, there is
no unified}\nonumber\\
&&\mbox{theory of partial differential equations.  Some equations
have their}\nonumber\\
&&\mbox{own theories, while others have no theory at all.  The
reason for}\nonumber\\
&&\mbox{this complexity is a more complicated
geometry...''}\nonumber
\end{eqnarray}
Within the confines of functional analysis, and more broadly
speaking, general topology, this was until very recently the case.
However, within the last fifteen years, two such general, type
independent theories for the solutions PDEs appeared.

In 1994 there appeared a type independent theory of existence and
regularity of solutions to nonlinear systems of PDEs based on
order completion of sets of functions, see \cite{Obergugenberger
and Rosinger}.  This theory yields the existence, and uniqueness,
of solutions to arbitrary continuous and nonlinear systems of
PDEs, with the solutions assimilated with Hausdorff continuous
functions, see \cite{AnguelovRosinger}.

Also, within the confines of functional analysis, there recently
emerged a general and type independent theory of PDEs, see
\cite{Neuberger 1} through \cite{Neuberger 4}. This theory,
developed by Neuberger, is based on approximations within Hilbert
space obtained by a generalized method of Steepest Descent.  This
method has yielded exceptional numerical results.

Both methods are in fact able to solve equations that are far more
general than PDEs.  This is precisely the feature which makes them
so powerful and type independent when applied to the particular
case of PDEs.

\subsection{A Brief Overview of the Order Completion Method}

Let us recall the general construction for solving a PDE of the
form
\begin{equation}\label{PDE1}
T\left(x,D\right)u\left(x\right)=f\left(x\right)\mbox{,
}x\in\Omega
\end{equation}
 by topological methods.  Here $\Omega\subseteq\mathbb{R}^{n}$ is supposed to be nonempty
and open.  The right hand term $f$ is assumed continuous, and the
unknown function satisfies, initially at least,
$u\in\mathcal{C}^{m}\left(\Omega\right)$.  The operator
$T\left(x,D\right)$ is supposed to be defined by some continuous
function $F:\Omega\times\mathbb{R}^{M}\rightarrow \mathbb{R}$
through
\begin{equation}\label{PDE2}
T\left(x,D\right)u\left(x\right)=F\left(x,u\left(x\right),...,D^{\alpha}u\left(x\right),...\right)\mbox{,
}x\in\Omega\mbox{ and }|\alpha|\leq m
\end{equation}
The general construction for finding generalized solutions to
(\ref{PDE1}) through (\ref{PDE2}), within the context of topology,
and specifically functional analysis, is summarized as follows:
\begin{eqnarray}
&\mbox{(1)  }&\mbox{Start with a PDE-operator $T$ mapping some
initial space of clas-}\nonumber\\
&&\mbox{sical functions $X$ into some space of functions $Y$ ,
with the right-}\nonumber\\
&&\mbox{hand term $f\in Y$.}\nonumber\\
&\mbox{(2)  }&\mbox{Define some structure, e.g. a uniformity, on
$Y$.}\nonumber\\
&\mbox{(3)  }&\mbox{Define a structure of the same sort on $X$
through pullback of the}\nonumber\\
&&\mbox{operator $T$ to obtain a space $X_{T}$ so that $T$ is
compatible with}\nonumber\\
&&\mbox{the customary structures on $X$ and $Y$.}\nonumber\\
&\mbox{(4)  }&\mbox{Construct the completions of $X_{T}$ and $Y$,
say
$X_{T}^{\sharp}$ and $Y^{\sharp}$, and extend}\nonumber\\
&&\mbox{$T$ to a mapping $T^{\sharp}:X_{T}^{\sharp}\rightarrow
Y^{\sharp}$}\nonumber\\
&\mbox{(5)  }&\mbox{If $f\in
T^{\sharp}\left(X_{T}^{\sharp}\right)$, one obtains a generalized
solution to the equation}\nonumber
\end{eqnarray}

The order completion method operates in a similar fashion, see
\cite{Obergugenberger and Rosinger}.  It is based on the
fundamental approximation result that, under a condition which is
necessary for the existence of a classical solution at
$x\in\Omega$, every continuous righthand term $f$ in (\ref{PDE1})
can be approximated from below (or from above) by functions in
$\mathcal{C}_{nd}^{m}\left(\Omega\right)$ where
\begin{equation}\label{CNDDef}
f:\Omega\rightarrow\mathbb{R}\in\mathcal{C}_{nd}^{m}\left(\Omega\right)\Leftrightarrow\left(\begin{array}{l}
  \exists\mbox{ }\Gamma_{f}\subset\Omega\mbox{ closed nowhere dense :} \\
  f\in\mathcal{C}^{m}\left(\Omega\setminus\Gamma_{f}\right) \\
\end{array}\right)
\end{equation}
In this regard, for any $f\in\mathcal{C}^{0}\left(\Omega\right)$
we have
\begin{eqnarray}\label{Approx1}
&&\forall\mbox{ }\epsilon
> 0\mbox{ :}\nonumber\\
&&\exists\mbox{ }\Gamma_{\epsilon}\subset\Omega\mbox{ closed
nowhere
dense, }U_{\epsilon}\in\mathcal{C}^{m}\left(\Omega\setminus\Gamma_{\epsilon}\right)\mbox{ :}\label{OCApprox}\\
&& f\left(x\right)-\epsilon\leq
T\left(x,D\right)U_{\epsilon}\left(x\right)\leq
f\left(x\right),x\in\Omega\setminus\Gamma_{\epsilon}\nonumber
\end{eqnarray}
Essentially, what this means is that the image of
$\mathcal{C}_{nd}^{m}\left(\Omega\right)$ under
$T\left(x,D\right)$ is order dense in
$\mathcal{C}_{nd}^{0}\left(\Omega\right)$.

One now proceeds as follows:  On the set
$\mathcal{C}^{0}_{nd}\left(\Omega\right)$ one considers the
equivalence relation
\begin{equation}\label{CNDEq}
f\sim g\Leftrightarrow\left(\begin{array}{l}
  \exists\mbox{ }\Gamma_{f,g}\subset\Omega\mbox{ closed, nowhere
dense}\mbox{ :} \\
  x\in\Omega\setminus\Gamma_{f,g}\Rightarrow f\left(x\right)=g\left(x\right) \\
\end{array}\right)
\end{equation}
On the quotient space
$\mathcal{M}^{0}\left(\Omega\right)=\mathcal{C}^{0}_{nd}\left(\Omega\right)/\sim$,
one introduces the partial order
\begin{equation}\label{CNDOrder}
F\leq G\Leftrightarrow\left(%
\begin{array}{l}
   \exists\mbox{ }f\in F\mbox{, }g\in G\mbox{, }\Gamma\subset\Omega\mbox{, closed, nowhere dense}\mbox{ :} \\
   (1)\mbox{  }f,g\in\mathcal{C}^{0}\left(\Omega\setminus\Gamma\right) \\
   (2)\mbox{  }f\leq g\mbox{ on }\Omega\setminus\Gamma \\
\end{array}%
\right)
\end{equation}
By the continuity of the operator $T\left(x,D\right)$ it follows
that if $u\in\mathcal{C}^{m}_{nd}\left(\Omega\right)$, and
$u\in\mathcal{C}^{m}\left(\Omega\setminus\Gamma\right)$, then
$T\left(x,D\right)u\in\mathcal{C}^{0}\left(\Omega\setminus\Gamma\right)$
so that one has the well defined mapping
\begin{equation}\label{ExtPDE1}
T\left(x,D\right):\mathcal{C}^{m}_{nd}\left(\Omega\right)\rightarrow\mathcal{C}^{0}_{nd}\left(\Omega\right)
\end{equation}

On the set $\mathcal{C}^{m}_{nd}\left(\Omega\right)$ one considers
the equivalence relation induced by the PDE operator
$T\left(x,D\right)$
\begin{eqnarray}\nonumber
u\sim_{T}v\Leftrightarrow T\left(x,D\right)u\sim
T\left(x,D\right)v.\nonumber
\end{eqnarray}
With the mapping (\ref{ExtPDE1}) one can now associate in a
canonical way an injective mapping
\begin{eqnarray}\nonumber
T:\mathcal{M}^{m}_{T}\left(\Omega\right)\rightarrow\mathcal{M}^{0}\left(\Omega\right),\nonumber
\end{eqnarray}
where
$\mathcal{M}^{m}_{T}\left(\Omega\right)=\mathcal{C}^{m}_{nd}\left(\Omega\right)/\sim_{T}$.
 Next one introduces a partial order on
$\mathcal{M}^{m}_{T}\left(\Omega\right)$ by
\begin{equation}\label{PulbackOrder}
U\leq_{T}V\Leftrightarrow TU\leq TV\mbox{ in
}\mathcal{M}^{0}\left(\Omega\right).
\end{equation}
This construction results in an order isomorphic embedding
\begin{eqnarray}\nonumber
T:\mathcal{M}^{m}_{T}\left(\Omega\right)\hookrightarrow
\mathcal{M}^{0}\left(\Omega\right).\nonumber
\end{eqnarray}
By extension of $T$ to the Dedekind completion
$\mathcal{M}_{T}^{m}\left(\Omega\right)^{\sharp}$ one obtains the
following commutative diagram

\begin{math}
\setlength{\unitlength}{1cm} \thicklines
\begin{picture}(13,6)
\put(2.4,5){$\mathcal{M}_{T}^{m}\left(\Omega\right)$}
\put(3.8,5.1){\vector(1,0){6.5}}
\put(10.5,5){$\mathcal{M}^{0}\left(\Omega\right)$}
\put(6.6,5.4){$T$} \put(3.0,4.8){\vector(0,-1){3.3}}
\put(2.4,1){$\mathcal{M}_{T}^{m}\left(\Omega\right)^{\sharp}$}
\put(4.0,1.1){\vector(1,0){6.2}} \put(6.6,1.3){$T^{\sharp}$}
\put(10.5,1){$\mathcal{M}^{0}\left(\Omega\right)^{\sharp}$}
\put(11.2,4.8){\vector(0,-1){3.3}}
\end{picture}
\end{math}
\\
\\
with $T^{\sharp}$ an order isomorphism.  Hence, for every
$f\in\mathcal{C}^{0}\left(\Omega\right)$, the generalized equation
\begin{eqnarray}
T^{\sharp}U=f\nonumber
\end{eqnarray}
has a unique solution in
$\mathcal{M}^{m}_{T}\left(\Omega\right)^{\sharp}$.

\subsection{Extended Concepts of Topology}

The aim of this paper is to cast the Order Completion Method of
Oberguggenberger and Rosinger within the framework of the
traditional topological and function analytical method.  In order
to achieve this aim it is necessary to consider topological type
structures that are more general than the usual
Hausdorff-Kuratowski-Bourbaki concept of topology.  Note that this
need for generalized concepts of topology is not unique to the
problem considered here, but appears frequently in analysis- even
in the relatively simple setting of locally convex linear spaces.
Indeed, recall that there is no topology on the topological dual
$E^{*}$ of a locally convex space $E$ such that the simple
evaluation mapping
\begin{eqnarray}
ev:E\times E^{*}\rightarrow\mathbb{R}\nonumber
\end{eqnarray}
defined through
\begin{eqnarray}
ev:\left(x,x^{*}\right)\mapsto\langle x^{*},x\rangle\nonumber
\end{eqnarray}
is jointly continuous, unless $E$ is a normable space.  Other
examples include some natural notions of convergence in measure
theory and in the theory of ordered spaces, which can not be
induced by a topology. For a more detailed motivation of the need
for generalized concepts of topology we refer the reader to
\cite{Beattie} and \cite{RosingervdWalt}.

The most general such topological type structures so far
considered in the literature was introduced by Rosinger, see
\cite{Rosinger2} through \cite{Rosinger4}, and more recently
\cite{RosingervdWalt}. A highly particular case of these
topological type structures are the celebrated convergence
structures and convergence spaces, see \cite{Beattie}, which are
defined in terms of convergent filters.  The collection of
convergence structures on a given set $X$ contains the collection
of all topologies on $X$ so that Hausdorff-Kuratowsky-Bourbaki
topology is a special case of convergence structures.

The theory of uniform spaces is a particular but highly important
case of general topology, among others due to the issue of
completeness and completion.  The generalization of uniform spaces
and uniformities within the setting of convergence spaces are
uniform convergence structures and uniform convergence spaces.
Whereas the uniform spaces are rather particular amongst
topological spaces, almost all convergence spaces are uniform
convergence space. Indeed, every Hausdorff convergence space is a
uniform convergence space, and there are even some uniform
convergence spaces that are not Hausdorff!  As is the case for
uniform spaces, every Hausdorff uniform convergence space admits a
Hausdorff completion which is unique up to homeomorphism.  For
details on uniform convergence spaces we refer the reader to
\cite{Beattie}, \cite{Gahler 1}, \cite{Gahler 2} or Appendix A.

\subsection{Order Convergence and the Order Convergence Structure}

As mentioned, there are many natural and useful modes of
convergence that can be defined in terms of a partial order $\leq$
on a set $X$, see for instance \cite{Birkhoff}, \cite{RieszI} and
\cite{Peressini}.  One of the better known such modes of
convergence is order convergence of sequences, which can be
defined on an arbitrary poset $X$ through
\begin{eqnarray}
\left(x_{n}\right)\mbox{ order converges to }x\in X
\Leftrightarrow\left(\begin{array}{l}
  \exists\mbox{ }\left(\lambda_{n}\right)\mbox{,
  }\left(\mu_{n}\right)\subset X\mbox{ :}\\
  \left(\lambda_{n}\right)\mbox{ increases to }x\mbox{ :} \\
  \left(\mu_{n}\right)\mbox{ decreases to }x\mbox{ :} \\
  n\in\mathbb{N}\Rightarrow \lambda_{n}\leq x_{n} \leq\mu_{n} \\
\end{array}\right)\label{OrderConvDef}
\end{eqnarray}
It is well know that order convergence of sequences can, in
general, not be induced by a topology, that is, for an arbitrary
poset $X$ there is no topology $\tau$ on $X$ such that the
convergent sequences with respect to $\tau$ are exactly those
sequences that order converge, see for instance \cite{AngvdWalt}.
However, recently it was shown that for a $\sigma$-distributive
lattice $X$ there is a convergence structure $\lambda$ on $X$ that
induces order convergence of sequences, see \cite{AngvdWalt} and
\cite{vdWalt1}.  One such a convergence structure, the Order
Convergence Structure $\lambda_{o}$, was defined as
\begin{eqnarray}
\mathcal{F}\in\lambda_{o}\left(x\right)\Leftrightarrow
\left(\begin{array}{l}
  \exists\mbox{ }\left(\lambda_{n}\right)\mbox{, } \left(\mu_{n}\right)\mbox{ on }X\mbox{ :}\\
  \left(\lambda_{n}\right)\mbox{ increases to }x\mbox{ :} \\
  \left(\mu_{n}\right)\mbox{ decreases to }x\mbox{ :} \\
  \{[\lambda_{n},\mu_{n}]\mbox{ : }n\in\mathbb{N}\} \subseteq\mathcal{F} \\
\end{array}\right)\label{OConStrucDef}
\end{eqnarray}
Note that this result applies amongst others to all vector
lattices.  Moreover, if a vector lattice $X$ is also Archimedean,
then the Order Convergence Structure is a vector space convergence
structure in the sense of \cite{Beattie}.  An important result
obtained in \cite{vdWalt1} is that the completion of the uniform
convergence structure induced on $X$ by $\lambda_{o}$ and the
linear structure is the Dedekind $\sigma$-completion of $X$
equipped with the Order Convergence Structure.  In case $X$ is
order separable \cite{RieszI}, the completion is in fact the
Dedekind completion of $X$.

Since the space $\mathcal{M}^{0}\left(\Omega\right)$ that appears
in the Order Completion Method can be viewed as an order separable
Archimedean vector lattice, the completion result described above
appears to provide a way of casting the Order Completion Method in
the framework of uniform convergence spaces.  However, one of the
main advantages of the Order Completion Method is that it makes no
distinction between linear and nonlinear equations, as is usually
the case for topological methods.  The completion result for the
Order Convergence Structure on an Archimedean vector lattice, on
the other hand, is a distinctly linear result.   However, in this
regard it was recently shown \cite{vdWalt2} that the Order
Convergence Structure on the order isomorphic copy
$\mathcal{ML}^{0}\left(\Omega\right)$ of
$\mathcal{M}^{0}\left(\Omega\right)$ is induced by a nontrivial
uniform convergence structure, and the resulting completion is
again the Dedekind completion $\mathbb{H}_{nf}\left(\Omega\right)$
of $\mathcal{ML}^{0}\left(\Omega\right)$.  Here
$\mathbb{H}_{nf}\left(\Omega\right)$ denotes the set of all nearly
finite Hausdorff continuous functions, \cite{Anguelov} and
\cite{Sendov}.

\section{Spaces of Functions}

This section is devoted to a short introduction of the spaces of
functions, and spaces of generalized functions on which the PDEs,
and their extensions, are to act.  As with the Order Completion
Method, we start with the space
$\mathcal{C}^{m}_{nd}\left(\Omega\right)$.  On this space we can
consider the Lower- and Upper Baire Operators
\begin{eqnarray}\label{IDef}
I\left(f\right)\left(x\right)=\sup\{\inf\{z\in
f\left(y\right)\mbox{ : }y\in V\}\mbox{ :
}V\in\mathcal{V}_{x}\}\mbox{, }x\in X
\end{eqnarray}
\begin{eqnarray}\label{SDef}
S\left(f\right)\left(x\right)=\inf\{\sup\{z\in
f\left(y\right)\mbox{ : }y\in V\}\mbox{ :
}V\in\mathcal{V}_{x}\}\mbox{, }x\in X
\end{eqnarray}
where $\mathcal{V}_{x}=\{V\subseteq X|V\mbox{ a neighbourhood of
}x\}$ is the neighborhood filter at $x\in\Omega$.  These operators
were first introduced by Baire \cite{Baire} for real valued
functions of a real variable, and later generalized by Sendov
\cite{Sendov} and Anguelov \cite{Anguelov}.  Note that the natural
setting of these operators, within the context of point valued
functions, is the set $\mathcal{A}\left(\Omega\right)$ consisting
of all functions $f$ that act on $\Omega$ and take values in the
extended real line
$\overline{\mathbb{R}}=\mathbb{R}\cup\{\pm\infty\}$.  There the
operators $I$ and $S$ characterize semi-continuity.  Indeed,
\begin{eqnarray}\label{LSCChar}
f\in\mathcal{A}\left(\Omega\right)\mbox{ is lower semi-continuous
}\Leftrightarrow I\left(f\right)=f
\end{eqnarray}
\begin{eqnarray}\label{USCChar}
f\in\mathcal{A}\left(\Omega\right)\mbox{ is upper semi-continuous
}\Leftrightarrow S\left(f\right)=f
\end{eqnarray}
Moreover, compositions of these operators characterize the normal
semi-continuous functions, see \cite{Dilworth} and
\cite{Anguelov}:
\begin{eqnarray}\label{NLSCChar}
f\in\mathcal{A}\left(\Omega\right)\mbox{ is normal lower
semi-continuous }\Leftrightarrow I\left(S\left(f\right)\right)=f
\end{eqnarray}
\begin{eqnarray}\label{NUSCChar}
f\in\mathcal{A}\left(\Omega\right)\mbox{ is normal upper
semi-continuous }\Leftrightarrow
S\left(I\left(f\right)\right)\left(x\right)=f
\end{eqnarray}
The normal lower and upper semi-continuous functions were
introduced in \cite{Dilworth} in connection with his attempts to
find the Dedekind order completion of spaces of continuous
functions.  This problem was finally solved by Anguelov, see
\cite{Anguelov}, where he characterized the Dedekind completion of
$\mathcal{C}\left(X\right)$ and some of its subspaces as spaces of
Hausdorff continuous functions.

We now consider the space
\begin{eqnarray}
\mathcal{ML}^{m}\left(\Omega\right)=\{I\left(S\left(f\right)\right)\mbox{
: }f\in\mathcal{C}^{m}_{nd}\left(\Omega\right)\}.\label{MLDef}
\end{eqnarray}
That is, $\mathcal{ML}^{m}\left(\Omega\right)$ is the image of
$\mathcal{C}^{m}_{nd}\left(\Omega\right)$ under the operator
$I\circ S$.  Since $\left(I\circ
S\right)\left(f\right)\left(x\right)$ is equal to
$f\left(x\right)$ whenever $f$ is continuous at $x$, it is $m$
times continuously differentiable everywhere except on a closed
nowhere dense subset of $\Omega$. Hence it follows that
$\mathcal{ML}^{m}\left(\Omega\right)\subset
\mathcal{C}^{m}_{nd}\left(\Omega\right)$.  By the idempotence of
the composite operator $I\circ S$, \cite{Anguelov}, it follows by
(\ref{NLSCChar}) and (\ref{CNDDef}) that
$\mathcal{ML}^{m}\left(\Omega\right)$ consists of all normal lower
semi-continuous functions which are $m$-times continuously
differentiable everywhere except on some closed nowhere dense
subset of $\Omega$. This space was introduced and studied in
\cite{vdWalt2} and \cite{vdWalt3}, and it was shown that
$\mathcal{ML}^{m}\left(\Omega\right)$ is an order isomorphic copy
of $\mathcal{M}^{m}\left(\Omega\right)$. Recall that
$\mathcal{M}^{m}\left(\Omega\right)=
\mathcal{C}^{m}_{nd}\left(\Omega\right)/\sim$, where $\sim$ is the
equivalence relation (\ref{CNDEq}).

We now consider the family $\mathcal{J}_{o}$ of filters on the
Cartesian product $\mathcal{ML}^{0}\left(\Omega\right)
\times\mathcal{ML}^{0}\left(\Omega\right)$.
\begin{definition}
A filter $\mathcal{U}$ on
$\mathcal{ML}^{0}\left(\Omega\right)\times\mathcal{ML}^{0}\left(\Omega\right)$
belongs to the family $\mathcal{J}_{o}$ whenever there exists
$k\in\mathbb{N}$ such that
\begin{eqnarray}\label{UOC}
&&\forall\mbox{ }i=1,...,k\mbox{ :}\nonumber\\
&&\exists\mbox{ }\left(I_{n}^{i}\right)\mbox{ nonempty order
intervals :}\nonumber\\
&&1)\mbox{  }I_{n}^{i}\supseteq I_{n+1}^{i}\mbox{ for every
}n\in\mathbb{N}\nonumber\\
&&2)\mbox{  }\forall\mbox{ }V\subseteq\Omega\mbox{ open
}\cap_{n\in\mathbb{N}}I_{n|V}^{i}=\emptyset\mbox{ or
there}\nonumber\\
&&\mbox{is }f_{i}\in \mathcal{ML}^{0}\left(V\right)\mbox{ such
that }\cap_{n\in\mathbb{N}}I_{n|V}^{i}=\{f_{i}\}\nonumber
\end{eqnarray}
and
\begin{eqnarray}
\left(\mathcal{I}^{1}\times\mathcal{I}^{1}\right)\cap...
\cap\left(\mathcal{I}^{k}\times\mathcal{I}^{k}\right)
\subset\mathcal{U}\mbox{ where }\mathcal{I}^{i}=[\{I^{i}_{n}\mbox{
: }n\in\mathbb{N}\}]\nonumber
\end{eqnarray}
\end{definition}
The family $\mathcal{J}_{o}$ is in fact a uniform convergence
structure on $\mathcal{ML}^{0}\left(\Omega\right)$, which we call
the Uniform Order Convergence Structure.  It induces the Order
Convergence Structure on $\mathcal{ML}^{0}\left(\Omega\right)$,
see \cite{vdWalt2} and \cite{vdWalt3}, and hence the following
properties are immediate:
\begin{eqnarray}
&\bullet&\mbox{A filter $\mathcal{F}$ on
$\mathcal{ML}^{0}\left(\Omega\right)$ converges with respect to
the induced convergence}\nonumber\\
&&\mbox{structure if and only if $\mathcal{F}$ converges in the
Order Convergence Structure}\nonumber\\
&&\mbox{on
$\mathcal{ML}^{0}\left(\Omega\right)$.}\label{UOCSP1}\\
&\bullet&\mbox{A sequence $\left(f_{n}\right)$ on
$\mathcal{ML}^{0}\left(\Omega\right)$ converges to
$f\in\mathcal{ML}^{0}\left(\Omega\right)$ if and only if
$\left(f_{n}\right)$}\nonumber\\
&&\mbox{order converges to $f$.}\label{UOCSP2}\\
&\bullet&\mbox{The uniform convergence space
$\mathcal{ML}^{0}\left(\Omega\right)$ is uniformly
Hausdorff.}\label{UOCSP3}
\end{eqnarray}

We now consider a system of $K$ nonlinear PDEs, each of order at
most $m$, in $K$ unknown functions $u_{1},...,u_{K}$
\begin{equation}\label{PDEDef}
\textbf{T}\left(x,D\right)\textbf{u}\left(x\right)=\textbf{f}\left(x\right)\mbox{,
}x\in\Omega\subseteq\mathbb{R}^{n}\mbox{ an open set.}
\end{equation}
Here $\textbf{u}:\Omega\rightarrow\mathbb{R}^{K}$ is a
$K$-dimensional vector valued function with components
$u_{1},...,u_{K}\in\mathcal{C}^{m}\left(\Omega\right)$.  The
righthand term $\textbf{f}:\Omega\rightarrow\mathbb{R}^{K}$ is
supposed to have continuous components
$f_{1},...,f_{K}:\Omega\rightarrow\mathbb{R}$.  The PDE-operator
$\textbf{T}\left(x,D\right)$ is defined through
\begin{equation}\label{PDEOperator}
\textbf{T}\left(x,D\right):\textbf{u}\left(x\right)\mapsto\textbf{F}\left(x,
...,u_{i}\left(x\right),...,D^{\alpha}u_{i},...\right)\mbox{ with
}i=1,...,K\mbox{ ,}|\alpha|\leq m
\end{equation}
where
$\textbf{F}:\Omega\times\mathbb{R}^{M}\rightarrow\mathbb{R}^{K}$
is a jointly continuous function in all of its arguments with
components
$F_{1},...,F_{K}:\Omega\times\mathbb{R}^{M}\rightarrow\mathbb{R}$.
One can now write equation (\ref{PDEDef}) in the form
\begin{eqnarray}\nonumber
\begin{array}{ccc}
  T_{1}\left(x,D\right)\textbf{u}\left(x\right) & = & f_{1}\left(x\right) \\
  \vdots & \vdots & \vdots \\
  T_{K}\left(x,D\right)\textbf{u}\left(x\right) & = & f_{k}\left(x\right) \\
\end{array}
\mbox{, }x\in\Omega.\nonumber
\end{eqnarray}
where each $T_{i}\left(x,D\right)$ is defined through
\begin{eqnarray}\nonumber
T_{i}\left(x,D\right):\textbf{u}\left(x\right)\mapsto
F_{i}\left(x,...,u_{j}\left(x\right),...,D^{\alpha}u_{j},...\right)\mbox{
with }j=1,...,K\nonumber
\end{eqnarray}

With the operator $T\left(x,D\right)$ we associate a mapping
\begin{eqnarray}
\widetilde{\textbf{T}}:\mathcal{ML}^{m}\left(\Omega\right)^{K}\rightarrow
\mathcal{ML}^{0}\left(\Omega\right)^{K}\nonumber
\end{eqnarray}
with components
\begin{eqnarray}
\widetilde{T}_{i}:\mathcal{ML}^{m}\left(\Omega\right)\rightarrow
\mathcal{ML}^{0}\left(\Omega\right)\nonumber
\end{eqnarray}
defined through
\begin{eqnarray}
\widetilde{T}_{i}\left(\textbf{u}\right):x\rightarrow \left(I\circ
S\right)\left(F\left(\cdot,u_{1},...,u_{K},...,\mathcal{D}^
{\alpha}u_{i},...\right)\right)\left(x\right)\mbox{,
}x\in\Omega\nonumber
\end{eqnarray}
where
\begin{eqnarray}
\mathcal{D}^{\alpha}f=\left(I\circ
S\right)\left(D^{\alpha}f\right)\mbox{,
}f\in\mathcal{C}^{m}_{nd}\left(\Omega\right).\nonumber
\end{eqnarray}
Since each $\mathcal{D}^{\alpha}u_{i}$ is continuous on an open
and dense set of $\Omega$, it follows by the continuity of the
$F_{i}$ that $\widetilde{\textbf{T}}\textbf{u}\in
\mathcal{ML}^{0}\left(\Omega\right)^{K}$ for every $\textbf{u}\in
\mathcal{ML}^{m}\left(\Omega\right)^{K}$.  Moreover, since
(\ref{LSCChar}) and (\ref{USCChar}) imply that a function
continuous at a point $x\in\Omega$ is invariant under the
operators $I$ and $S$, at $x$, it follows that
\begin{eqnarray}
&&\forall\mbox{
}\textbf{u}\in\mathcal{ML}^{m}\left(\Omega\right)\mbox{
:}\nonumber\\
&&\exists\mbox{ }\Gamma\subset\Omega\mbox{ :}\label{TExtTRel}\\
&&\Gamma\mbox{ is closed nowhere dense :}\nonumber\\
&&x\in\Omega\setminus\Gamma\Rightarrow
\widetilde{\textbf{T}}\left(\textbf{u}\right)\left(x\right)=
\textbf{T}\left(x,D\right)\textbf{u}\left(x\right)\nonumber
\end{eqnarray}
On the space $\mathcal{ML}^{m}\left(\Omega\right)^{K}$ consider
the equivalence relation induced by $\widetilde{\textbf{T}}$
through
\begin{equation}\label{TEquiv}
\textbf{u}\sim_{\widetilde{\textbf{T}}}\textbf{v}\Leftrightarrow
\widetilde{\textbf{T}}\textbf{u}=\widetilde{\textbf{T}}\textbf{v}
\end{equation}
We denote the quotient space
$\mathcal{ML}^{m}\left(\Omega\right)^{K}/\sim_{\widetilde{\textbf{T}}}$
by $\mathcal{ML}^{m}_{\widetilde{\textbf{T}}}\left(\Omega\right)$.
With the mapping $\widetilde{\textbf{T}}$ we now associate a
mapping
\begin{equation}\label{PulbackOp}
\widehat{\textbf{T}}:\mathcal{ML}^{m}_{\widetilde{\textbf{T}}}\left(\Omega\right)
\rightarrow \mathcal{ML}^{0}\left(\Omega\right)^{K}
\end{equation}
which assigns to each equivalence class in
$\mathcal{ML}^{m}_{\widetilde{\textbf{T}}}\left(\Omega\right)$ its
image under the operator $\widetilde{\textbf{T}}$, that is,
\begin{eqnarray}
&&\forall\mbox{ }\textbf{U}\in\mathcal{ML}^{m}
_{\widetilde{\textbf{T}}}\left(\Omega\right)\mbox{ :}\nonumber\\
&&\forall\mbox{ }\textbf{u}\in\textbf{U}\mbox{ :}\label{TExtTRel2}\\
&&\widehat{\textbf{T}}\left(\textbf{U}\right)=\widetilde{\textbf{T}}\left(\textbf{u}\right)\nonumber
\end{eqnarray}
The generalized version of equation (\ref{PDEDef}) is
\begin{equation}\label{GenPDE}
\widehat{\textbf{T}}\textbf{U}=\textbf{f},
\end{equation}
where $\textbf{U}$ is an equivalence class under the equivalence
$\sim_{\widetilde{\textbf{T}}}$.  The mapping
$\widehat{\textbf{T}}$ is injective, but it is in general not
surjective.

We consider the space $\mathcal{ML}^{0}\left(\Omega\right)^{K}$
equipped with the product uniform convergence structure when each
copy of $\mathcal{ML}^{0}\left(\Omega\right)$ carries the Uniform
Order Convergence Structure.
\begin{definition}\label{ProdUCSML0}
A filter $\mathcal{U}$ on
$\mathcal{ML}^{0}\left(\Omega\right)^{K}\times
\mathcal{ML}^{0}\left(\Omega\right)^{K}$ belongs to
$\mathcal{J}^{K}_{o}$ whenever $\left(\pi_{i}\times\pi_{i}\right)
\left(\mathcal{U}\right)\in\mathcal{J}_{o}$ for every $i=1,...,K$.
Here $\pi_{i}:\mathcal{ML}^{0}\left(\Omega\right)^{K} \rightarrow
\mathcal{ML}^{0}\left(\Omega\right)$ denotes the projection along
the $i$th coordinate.
\end{definition}
$\mathcal{ML}^{m}_{\widetilde{\textbf{T}}}\left(\Omega\right)$
carries the initial uniform convergence structure
$\mathcal{J}_{\widetilde{\textbf{T}}}$ with respect to the
injective mapping
\begin{eqnarray}\label{InjMap}
\widehat{\textbf{T}}:\mathcal{ML}^{m}_{\widetilde{\textbf{T}}}\left(\Omega\right)
\rightarrow \mathcal{ML}^{0}\left(\Omega\right)^{K}
\end{eqnarray}
\begin{definition}\label{InitCSMLmT}
A filter $\mathcal{U}$ on $\mathcal{ML}^{m}_
{\widetilde{\textbf{T}}}\left(\Omega\right)\times
\mathcal{ML}^{m}_{\widetilde{\textbf{T}}}\left(\Omega\right)$
belongs to $\mathcal{J}_{\widetilde{\textbf{T}}}$ whenever
$\left(\widetilde{\textbf{T}}\times
\widetilde{\textbf{T}}\right)\left(\mathcal{U}\right)$ belongs to
$\mathcal{J}^{K}_{o}$.
\end{definition}
From the injectivity of $\widehat{\textbf{T}}$  we obtain
\begin{eqnarray}
&&\widehat{\textbf{T}}:
\mathcal{ML}^{m}_{\widetilde{\textbf{T}}}\left(\Omega\right)
\hookrightarrow \mathcal{ML}^{0}\left(\Omega\right)^{K}\mbox{ is a
uniformly continuous embedding}\label{TEmbedding}
\end{eqnarray}

The completion of $\mathcal{ML}^{0}\left(\Omega\right)$ can be
constructed as a set of Hausdorff continuous functions.  In fact,
it is the set $\mathbb{H}_{nf}\left(\Omega\right)$ consisting of
all nearly finite Hausdorff continuous functions, equipped with a
suitably chosen uniform convergence structure, see \cite{vdWalt3}
for details.  Hence it follows by Theorem \ref{ProdCom} that the
completion of $\mathcal{ML}^{0}\left(\Omega\right)^{K}$ is
$\mathbb{H}_{nf}\left(\Omega\right)^{K}$, equipped with the
product uniform convergence structure.  Moreover, it now follows
by (\ref{TEmbedding}) and Theorem \ref{SubspUCSComp} that the
completion of $\mathcal{ML}^{m}_{\widetilde{\textbf{T}}}
\left(\Omega\right)$ is homeomorphic to a subspace
$\mathbb{H}_{\widetilde{\textbf{T}}}\left(\Omega\right)$ of
$\mathbb{H}_{nf}\left(\Omega\right)$.  We will not give details
concerning Hausdorff continuous functions.  For details concerning
Hausdorff continuous functions, we refer the reader to
\cite{Anguelov}.

\section{Approximation Results}

We now again consider a system of $K$ nonlinear PDEs of the form
(\ref{PDEDef}) through (\ref{PDEOperator}). Recall that the Order
Completion Method for single nonlinear PDEs of the form
(\ref{PDE1}) through (\ref{PDE2}) is based on the simple
approximation result (\ref{Approx1}).  This result is follows from
the fact that every continuous function can be approximated
locally by a sequence of polynomials.  In this section we extend
this result to the general $K$-dimensional case, for $K\geq 1$
arbitrary but given.

A natural assumption on the function
$\textbf{F}:\Omega\times\mathbb{R}^{M}\rightarrow\mathbb{R}^{K}$,
and hence the PDE-operator $\textbf{T}\left(x,D\right)$, and the
righthand term $\textbf{f}$ is that, for every $x\in\Omega$
\begin{equation}\label{Assumption}
\textbf{f}\left(x\right)\in\textrm{int}R_{x}\mbox{ where
}R_{x}=\{\textbf{F}\left(x,\xi_{1\alpha},...,\xi_{i\alpha},...\right)\mbox{
: }\xi_{i\alpha}\in\mathbb{R}\mbox{, }i=1,...,K\mbox{,
}|\alpha|\leq m\}
\end{equation}
Note that (\ref{Assumption}) is is of a technical nature, and
hardly a restriction on the class of PDEs considered.  In fact,
every linear PDE, and also most nonlinear PDEs, satisfy
(\ref{Assumption}).  It is necessary condition for the existence
of a classical solution to (\ref{PDEDef}) in a neighborhood of
$x$.
\begin{theorem}\label{Approx}
Consider a system of PDEs of the form (\ref{PDEDef}) through
(\ref{PDEOperator}) that also satisfies (\ref{Assumption}).  For
every $\epsilon>0$ there exists a closed nowhere dense set
$\Gamma_{\epsilon}\subset\Omega$ with zero Lebesgue measure, and a
function
$\textbf{U}_{\epsilon}\in\mathcal{C}^{m}\left(\Omega\setminus\Gamma_{\epsilon}\right)^{K}$
with components $U_{\epsilon,1},...,U_{\epsilon,K}$ such that
\begin{equation}\label{ApEq}
f_{i}\left(x\right)-\epsilon\leq
T_{i}\left(x,D\right)\textbf{U}_{\epsilon}\left(x\right)\leq
f_{i}\left(x\right)\mbox{, }x\in\Omega\setminus\Gamma_{\epsilon}
\end{equation}
\end{theorem}
\begin{proof}\\
Let
\begin{equation}\label{AE3}
\Omega=\bigcup_{\nu\in\mathbb{N}}C_{\nu}
\end{equation}
where, for $\nu\in\mathbb{N}$, the compact sets $C_{\nu}$ are
$n$-dimensional intervals
\begin{equation}
C_{\nu}=[a_{\nu},b_{\nu}]
\end{equation}
with $a_{\nu}=\left(a_{\nu,1},...,a_{\nu,n}\right)$,
$b_{\nu}=\left(b_{\nu,1},...,b_{\nu,n}\right)\in\mathbb{R}^{n}$
and $a_{\nu,i}\leq b_{\nu,i}$ for every $i=1,...,n$.We also assume
that $C_{\nu}$, with $\nu\in\mathbb{N}$ are locally finite, that
is,
\begin{equation}\label{AE2}
\forall x\in\Omega\exists V_{x}\subseteq\Omega\mbox{ a
neighborhood of }x\mbox{ : }\{\nu\in\mathbb{N}\mbox{ :
}C_{\nu}\cap V_{x}\neq\emptyset\}\mbox{ is finite}
\end{equation}
We also assume that the interiors of $C_{\nu}$, with
$\nu\in\mathbb{N}$, are pairwise disjoint.  We note that note that
such $C_{\nu}$ exist, see \cite{Forster}.\\
Let us now take $\epsilon>0$ given arbitrary but fixed.  Let us
take $\nu\in\mathbb{N}$ and apply Lemma \ref{ApproxLemma} to each
$x_{0}\in C_{\nu}$.  Then we obtain $\delta_{x_{0}}>0$ and
$P_{x_{0},1},...,P_{x_{0},K}$ polynomial in $x\in\mathbb{R}^{n}$
such that
\begin{equation}
f_{i}\left(x\right)-\epsilon\leq
T_{i}\left(x,D\right)\textbf{P}_{x_{0}}\left(x\right)\leq
f\left(x\right)\mbox{,
}x\in\Omega\cap\overline{B}\left(x_{0},\delta_{x_{0}}\right)\mbox{
and }i=1,...,K
\end{equation}
where $\textbf{P}_{x_{0}}:\mathbb{R}^{n}\rightarrow\mathbb{R}^{K}$
is the $K$-dimensional vector valued function with components
$P_{x_{0},1},...,P_{x_{0},K}$.  Since $C_{\nu}$ is compact, it
follows that
\begin{eqnarray}
&&\exists\delta>0\mbox{ : }\nonumber\\
&&\forall x_{0}\in C_{\nu}\mbox{ : }\nonumber\\
&&\exists P_{x_{0},1},...,P_{x_{0},K}\mbox{ polynomial in }x\in\mathbb{R}^{n}\mbox{ : }\label{AE}\\
&&\|x-x_{0}\|\leq\delta\Rightarrow
f_{i}\left(x\right)-\epsilon\leq
T_{i}\left(x,D\right)\textbf{P}_{x_{0}}\left(x\right)\leq
f\left(x\right)\mbox{,
}x\in\overline{B}\left(x_{0},\delta\right)\cap \Omega\nonumber
\end{eqnarray}
where $i=1,...,K$.  Subdivide $C_{\nu}$ into $n$-dimensional
intervals $I_{\nu,1},...,I_{\nu,\mu}$ with diameter not exceeding
$\delta$ such that their interiors are pairwise disjoint.  If
$a_{j}$ with $j=1,...,\mu$ is the center of the interval
$I_{\nu,j}$ then by (\ref{AE}) there exists
$P_{a_{j},1},...,P_{a_{j}},K$ polynomial in $x\in\mathbb{R}^{n}$
such that
\begin{equation}\label{AE1}
f_{i}\left(x\right)-\epsilon\leq
T_{i}\left(x,D\right)\textbf{P}_{a_{j}}\left(x\right)\leq
f_{i}\left(x\right)\mbox{, }x\in I_{\nu,j}
\end{equation}
where $i=1,...,K$.  Now set
\begin{equation}
\Gamma_{\nu,\epsilon}=C_{\nu}\setminus\left(\left(\bigcup_{j=1}
^{\mu}\textrm{int}I_{\nu,j}\right)\cup\textrm{int}C_{\nu}\right)
\end{equation}
that is, $\Gamma_{\nu,\epsilon}$ is a rectangular grid generated
as a finite union of hyperplanes.  Furthermore, using (\ref{AE1}),
we find
\begin{equation}\label{AE4}
\textbf{U}_{\nu,\epsilon}\in\mathcal{C}^{m}\left(C_{\nu}\setminus\Gamma_{\nu,\epsilon}\right)
\end{equation}
such that
\begin{equation}\label{AE5}
f_{i}\left(x\right)-\epsilon\leq
T_{i}\left(x,D\right)\textbf{U}_{\nu,\epsilon}\left(x\right)\leq
f_{i}\left(x\right)\mbox{, }x\in
C_{\nu}\setminus\Gamma_{\nu,\epsilon}
\end{equation}
In view of (\ref{AE2}) it follows that
\begin{equation}
\Gamma_{\epsilon}=\bigcup_{\nu\mathbb{N}}\Gamma_{\nu,\epsilon}\mbox{
is closed nowhere dense and
}\textrm{mes}\left(\Gamma_{\epsilon}\right)=0
\end{equation}
From (\ref{AE3}), (\ref{AE4}) and (\ref{AE5}) we obtain
(\ref{ApEq}).
\end{proof}\\
The above proof relies on the following lemma which is in fact the
basic approximation result.
\begin{lemma}\label{ApproxLemma}
Consider a system of PDEs of the form (\ref{PDEDef}) through
(\ref{PDEOperator}) that also satisfies (\ref{Assumption}).  Then
\begin{eqnarray}
&&\forall\mbox{ }x_{0}\in\Omega\mbox{, }\epsilon>0\mbox{
:}\nonumber\\
&&\exists\mbox{ }\delta>0\mbox{, }P_{1},...,P_{K}\mbox{ polynomial
in $x\in\mathbb{R}^{n}$ :}\label{ALemEq}\\
&&x\in B\left(x_{0},\delta\right)\cap\Omega\mbox{, }1\leq i\leq k
\Rightarrow f_{i}\left(x\right)-\epsilon\leq
T_{i}\left(x,D\right)\textbf{P}\left(x\right)\leq
f_{i}\left(x\right)\nonumber
\end{eqnarray}
Here $\textbf{P}$ is the $K$-dimensional vector valued function
with components $P_{1},...,P_{K}$.
\end{lemma}
\begin{proof}\\
For any $x_{0}\in\Omega$ and $\epsilon>0$ small enough it follows
by (\ref{Assumption}) that there exist
\begin{equation}
\xi_{i\alpha}\in\mathbb{R}\mbox{ with }i=1,...,K\mbox{ and
}|\alpha|\leq m
\end{equation}
such that
\begin{equation}
F_{i}\left(x_{0},...,\xi_{i\alpha},...\right)=f_{i}\left(x_{0}\right)-\frac{\epsilon}{2}
\end{equation}
Now take $P_{1},...,P_{K}$ polynomials in $x\in\mathbb{R}^{n}$
that satisfy
\begin{equation}
D^{\alpha}P_{i}\left(x_{0}\right)=\xi_{i\alpha}\mbox{ for
}i=1,...,K\mbox{ and }|\alpha|\leq m
\end{equation}
Then it is clear that
\begin{equation}
T_{i}\left(x,D\right)\textbf{P}\left(x_{0}\right)-f_{i}\left(x_{0}\right)=-\frac{\epsilon}{2}
\end{equation}
where $\textbf{P}$ is the $K$-dimensional vector valued function
on $\mathbb{R}^{n}$ with components $P_{1},...,P_{K}$.  Hence
(\ref{ALemEq}) follows by the continuity of the $f_{i}$ and the
$F_{i}$.
\end{proof}\\

\section{Existence of Generalized Solutions}

We now proceed to establish both the existence and uniqueness of
generalized solutions to the generalized equation (\ref{GenPDE}).
This result appears as an application of the theory developed in
Sections 2 and 3, and the results concerning completions of
uniform convergence spaces contained in Appendix A.

The uniqueness of the solution should not be interpreted as
implying that any meaningful solutions are disregarded.  Indeed,
recall that the elements of the space
$\mathcal{ML}^{m}_{\widetilde{\textbf{T}}}\left(\Omega\right)$ are
equivalence classes under the equivalence (\ref{TEquiv}).  Hence,
if there is, for instance, a solution $\textbf{U}\in
\mathcal{ML}^{m}_{\widetilde{\textbf{T}}}\left(\Omega\right)$,
then $\textbf{U}$ is the equivalence class containing all possible
solutions $\textbf{u}\in\mathcal{ML}^{m}\left(\Omega\right)^{K}$
to (\ref{PDEDef}).

Consider now the following commutative diagram:\\
\\
\\
\begin{math}
\setlength{\unitlength}{1cm} \thicklines
\begin{picture}(13,6)

\put(3.4,5.4){$\mathcal{ML}^{m}_{\widetilde{\textbf{T}}}\left(\Omega\right)$}
\put(5.1,5.5){\vector(1,0){6.4}}
\put(11.6,5.4){$\mathcal{ML}^{0}\left(\Omega\right)^{K}$}
\put(7.8,5.9){$\widehat{\textbf{T}}$}
\put(3.8,5.2){\vector(0,-1){3.5}} \put(4.7,1.2){\vector(1,0){6.7}}
\put(3.4,1.1){$\mathbb{H}_{\widehat{\textbf{T}}}\left(\Omega\right)$}
\put(11.6,1.1){$\mathbb{H}_{nf}\left(\Omega\right)^{K}$}
\put(3.5,3.4){$\phi$} \put(12.2,3.4){$\varphi$}
\put(12.0,5.2){\vector(0,-1){3.5}}
\put(7.8,1.4){$\widehat{\textbf{T}}^{\sharp}$}

\end{picture}
\end{math}\\
\\
Here $\phi$ and $\varphi$ are the uniformly continuous embeddings
associated with the completions
$\mathbb{H}_{\widehat{\textbf{T}}}\left(\Omega\right)$ and
$\mathbb{H}_{nf}\left(\Omega\right)$, and
$\widehat{\textbf{T}}^{\sharp}$ is the extension of
$\widehat{\textbf{T}}$ achieved through uniform continuity.  The
existence of generalized solutions will follow from Theorem
\ref{Approx1} and (\ref{TEmbedding}).

\begin{theorem}\label{SolEx}
For every $\textbf{f}\in\mathcal{C}^{0}\left(\Omega\right)^{K}$
that satisfies (\ref{Assumption}), there exists a unique
$\textbf{U}\in \mathbb{H}_{\widetilde{\textbf{T}}}
\left(\Omega\right)^{K}$ such that
\begin{equation}\label{ExtPDEII}
\widehat{\textbf{T}}^{\sharp}\textbf{U}=\textbf{f}
\end{equation}
\end{theorem}
\begin{proof}
First let us show existence.  For every $n\in\mathbb{N}$, Theorem
\ref{Approx} yields a closed nowhere dense set
$\Gamma_{n}\subset\Omega$ and a function
$\textbf{u}_{n}\in\mathcal{C}^{m}\left(\Omega\right)$ that
satisfies
\begin{equation}
x\in\Omega\setminus\Gamma_{n}\Rightarrow
f_{i}\left(x\right)-\frac{1}{n}\leq
T_{i}\left(x,D\right)\textbf{u}_{n}\left(x\right)\leq
f_{i}\left(x\right)\mbox{, }i=1,...,K\label{Approx2}
\end{equation}
Since $\Gamma_{n}$ is closed nowhere dense we associate
$\textbf{u}_{n}$ with a function
$\textbf{v}_{n}\in\mathcal{ML}^{m}\left(\Omega\right)$ in a unique
way.  Indeed, consider for instance the function
\begin{eqnarray}\nonumber
\textbf{w}_{n}:x\mapsto\left\{\begin{array}{ll}
  \textbf{u}_{n}\left(x\right) & \mbox{ if }x\in\Omega\setminus\Gamma \\
  0 & \mbox{ if }x\in\Gamma \\
\end{array}\right.
\end{eqnarray}
Now let $\textbf{v}_{n}$ be the $K$-dimensional vector valued
function with components $v^{i}_{n}=\left(I\circ
S\right)\left(w_{n}^{i}\right)$.\\
Denote by $\textbf{V}_{n}$ the
equivalence class generated by $\textbf{v}_{n}$ under the
equivalence relation $\sim_{\widetilde{\textbf{T}}}$.  Applying
(\ref{TExtTRel}), (\ref{TExtTRel2}) and (\ref{Approx2}) one finds
that the sequence
$\left(\widehat{\textbf{T}}\left(\textbf{V}_{n}\right)_{i}\right)$
is increasing in $\mathcal{ML}^{0}\left(\Omega\right)$ and has
upper bound $f_{i}$, for each $i=1,...,K$.  It is clear that for
each $i=1,...,K$, $f_{i}$ is in fact the least upper of the
sequence
$\left(\widehat{\textbf{T}}\left(\textbf{V}_{n}\right)_{i}\right)$.
Hence each $\left(\widehat{\textbf{T}}
\left(\textbf{V}_{n}\right)_{i}\right)$ converges to $f_{i}$ in
$\mathcal{ML}^{0}\left(\Omega\right)$ so that the sequence
$\left(\widehat{\textbf{T}}\left(\textbf{V}_{n}\right)\right)$
converges to $\textbf{f}$ in
$\mathcal{ML}^{0}\left(\Omega\right)^{K}$.  It now follows that
$\left(\textbf{V}_{n}\right)$ is a Cauchy sequence in
$\mathcal{ML}^{m}_{\widetilde{\textbf{T}}}\left(\Omega\right)$ so
that there exists
$\textbf{U}\in\mathbb{H}_{\widehat{\textbf{T}}}\left(\Omega\right)$
that satisfies (\ref{ExtPDEII}).\\
Since the mapping
$\widehat{\textbf{T}}:\mathcal{ML}^{m}_{\widetilde{\textbf{T}}}\left(\Omega\right)
\rightarrow\widehat{\textbf{T}}$ is a uniformly continuous
embedding by (\ref{TEmbedding}), the uniqueness of the solution
$\textbf{U}$ found above now follows by Theorem \ref{EmbExt}.
\end{proof}

\section{Conclusion}

It has been shown that the Order Completion Method for systems of
nonlinear PDEs can be naturally considered within the framework of
uniform convergence spaces.  With respect to the previous
existence and regularity results obtained by the Order Completion
Method, we have not made any improvement.  The class of systems of
nonlinear PDEs treated is the same as that considered in
\cite{Obergugenberger and Rosinger}, that is, we considered all
continuous PDEs.

However, the generalized solutions obtained are elements of the
completion of a space of equivalence classes of piecewise smooth
functions, which are assimilated with Hausdorff continuous
functions.  In order to obtain meaningful solutions, for instance
classical solutions, a great deal of hard analysis is required.
Within the framework of uniform convergence spaces as applied
here, we have a sound topological foundation at our disposal.

\appendix

\section{Convergence Structures}

A filter on a set $X$ is a nonempty collection $\mathcal{F}$ of
nonempty subsets of $X$ that is closed under the formation of
supersets and finite intersection.  For any $x\in X$
\begin{eqnarray}
[x]=\{F\subseteq X\mbox{ : }x\in F\}\nonumber
\end{eqnarray}
is the filter generated by $x$.  A convergence structure on $X$ is
now defined as follows.
\begin{definition}\label{CSDef}
A convergence structure on $X$ is a mapping $\lambda$ from $X$
into the powerset of the set of all filters on $X$ that satisfies
\begin{eqnarray}
&\mbox{(1) }&[x]\in\lambda\left(x\right)\mbox{ for every }x\in
X\nonumber\\
&\mbox{(2) }&\mbox{If }\mathcal{F}\mbox{, }\mathcal{G}
\in\lambda\left(x\right)\mbox{, then
}\mathcal{F}\cap\mathcal{G}\in\lambda\left(x\right)\nonumber\\
&\mbox{(3) }&\mbox{If }\mathcal{F}\in\lambda\left(x\right)\mbox{
and }\mathcal{F}\subseteq\mathcal{G}\mbox{ then
}\mathcal{G}\in\lambda\left(x\right)\nonumber
\end{eqnarray}
The pair $\left(X,\lambda\right)$ is called a convergence space.
\end{definition}
A convergence space $X=\left(X,\lambda\right)$ is called Hausdorff
if it satisfies
\begin{eqnarray}
x\neq y\Rightarrow\lambda\left(x\right)\cap\lambda\left(y\right)
=\emptyset\label{Hausdorff}
\end{eqnarray}
A convergence structure $\lambda$ on $X$ induces convergence of
sequences through
\begin{eqnarray}
\left(x_{n}\right)\mbox{ converges to }x\in X\Leftrightarrow
\langle\left(x_{n}\right)\rangle=[\{\{x_{n}\mbox{ : }n\geq
k\}\mbox{ : }k\in\mathbb{N}\}]\in
\lambda\left(x\right)\label{SeqConv}
\end{eqnarray}
One can construct new convergence structures from old ones in the
following way.
\begin{definition}\label{InitCS}
Let $\left(X_{i}\right)_{i\in I}$ be a family of convergence
spaces, $X$ a set and for every $i\in I$ let $f_{i}:X\rightarrow
X_{i}$ be a mapping.  A filter $\mathcal{F}$ on $X$ converges to
$x\in X$ with respect to the initial convergence structure
$\lambda$ on $X$ with respect to $\left(f_{i}\right)_{i\in I}$
whenever
\begin{eqnarray}
f_{i}\left(\mathcal{F}\right)\in\lambda_{i}\left(x\right)\mbox{
for every }i\in I\nonumber
\end{eqnarray}
\end{definition}
This construction leads, as a particular case, to the product and
subspace convergence structures.

In order to define a uniform convergence structure we recall some
notation.  For subsets $U$ and $V$ of $X\times X$
\begin{eqnarray}\nonumber
U^{-1}=\{\left(x,y\right)\in X\times X\mbox{ :
}\left(y,x\right)\in U\}\nonumber
\end{eqnarray}
and the composition of $U$ and $V$ is given by
\begin{eqnarray}\nonumber
U\circ V=\{\left(x,y\right)\in X\times X\mbox{ : }\begin{array}{l}
  \exists\mbox{ }z\in X\mbox{ :} \\
  \left(x,z\right)\in V,\left(z,y\right)\in
U \\
\end{array}\}\nonumber
\end{eqnarray}
If $F\subseteq X$ then
\begin{eqnarray}\nonumber
U[F]=\{x\in X\mbox{ : }\begin{array}{l}
  \exists\mbox{ }y\in X \\
  \left(y,x\right)\in U \\
\end{array}\}\nonumber
\end{eqnarray}
If $x\in X$, then one sets
\begin{eqnarray}\nonumber
U[x]=U[\{x\}]\nonumber
\end{eqnarray}
For filters $\mathcal{U}$ and $\mathcal{V}$ we set
\begin{eqnarray}\nonumber
\mathcal{U}^{-1}=\{U^{-1}\mbox{ : }U\in\mathcal{U}\}\nonumber
\end{eqnarray}
and
\begin{eqnarray}\nonumber
\mathcal{U}\circ\mathcal{V}=[\{U\circ V\mbox{ :
}U\in\mathcal{U},V\in\mathcal{V}\}]\nonumber
\end{eqnarray}
whenever $U\circ V\neq\emptyset$ for every $U\in\mathcal{U}$ and
$V\in\mathcal{V}$.  If $\mathcal{F}$ is a filter on $X$, then
\begin{eqnarray}\nonumber
\mathcal{U}[\mathcal{F}]=[\{U[F]\mbox{ :
}U\in\mathcal{U},F\in\mathcal{F}\}]\nonumber
\end{eqnarray}
provided that $U[F]\neq\emptyset$ for every $U\in\mathcal{U}$ and
$F\in\mathcal{F}$.

A uniform convergence structure on $X$ is a family $\mathcal{J}$
of filters on $X\times X$, as apposed to a uniformity on $X$ which
consists of a single filter on $X\times X$.
\begin{definition}\label{UCSDef}
Let $X$ be a set.  A family $\mathcal{J}$ of filters on $X\times
X$ is called a uniform convergence structure on $X$ and
$\left(X,\mathcal{J}\right)$ a uniform convergence space if the
following hold:
\begin{eqnarray}
&\mbox{(1) }&\mbox{$[x]\times [x]\in\mathcal{J}$ for all $x\in X$.}\nonumber\\
&\mbox{(2) }&\mbox{$\mathcal{U}\cap\mathcal{V}\in\mathcal{J}$
whenever
$\mathcal{U},\mathcal{V}\in\mathcal{J}$.}\nonumber\\
&\mbox{(3) }&\mbox{If $\mathcal{U}\in\mathcal{J}$ then
$\mathcal{V}\in\mathcal{J}$ for every filter $\mathcal{V}$ on
$X\times X$ such that $\mathcal{V}\supseteq\mathcal{U}$.}\nonumber\\
&\mbox{(4) }&\mbox{If $\mathcal{U}\in\mathcal{J}$ then
$\mathcal{U}^{-1}\in\mathcal{J}$.}\nonumber\\
&\mbox{(5) }&\mbox{For all $\mathcal{U},\mathcal{V}\in\mathcal{J}$
one has $\mathcal{U}\circ\mathcal{V}\in\mathcal{J}$ whenever the
composition $\mathcal{U}\circ\mathcal{V}$ exists.}\nonumber
\end{eqnarray}
\end{definition}
A uniform convergence structure $\mathcal{J}$ on $X$ induces a
convergence structure $\lambda_{\mathcal{J}}$, the induced
convergence structure, on $X$ through
\begin{eqnarray}
\mathcal{F}\in\lambda_{\mathcal{J}}\left(x\right)\Leftrightarrow
[x]\times\mathcal{F}\in\mathcal{J}\label{IndCSDef}
\end{eqnarray}
We can construct new uniform convergence structures from old ones
in the following way:
\begin{definition}\label{InitUCSDef}
Let $X$ be a set and $\left(X_{i}\right)_{i\in I}$ a family of
uniform convergence spaces, and for each $i\in I$ let
$f_{i}:X\rightarrow X_{i}$ be a mapping.  The initial uniform
convergence structure $\mathcal{J}$ on $X$ with respect to
$\left(f_{i}\right)_{i\in I}$ is defined through
\begin{eqnarray}
\mathcal{U}\in\mathcal{J}\Leftrightarrow \left(f_{i}\times
f_{i}\right)\left(\mathcal{U}\right)\in\mathcal{J}_{i}\nonumber
\end{eqnarray}
\end{definition}
This general construction yields amongst others the product and
subspace uniform convergence structures.  The product convergence
structure on the Cartesian product $X$ of uniform convergence
spaces $\left(X_{i}\right)_{i\in I}$ is the initial uniform
convergence structure with respect to the projection mappings
$\left(\pi_{i}:X\rightarrow X_{i}\right)$.  The subspace uniform
convergence structure on a subset $Y$ of a uniform convergence
structure $X$ is the initial uniform convergence structure with
respect to the inclusion mapping $\iota:Y\rightarrow X$.  The
following result \cite{Beattie} is useful.
\begin{proposition}\label{InitCSInitUCS}
Let $X$ carry the initial uniform convergence structure with
respect to $\left(f_{i}:X\rightarrow X_{i}\right)$.  Then the
induced convergence structure is the initial convergence structure
with respect to the induced convergence structures on the $X_{i}$.
\end{proposition}
Some important classes of mappings between uniform convergence
spaces are the following.
\begin{definition}\label{UConFunDef}
Let $f:X\rightarrow Y$ be a mapping from the uniform convergence
space $\left(X,\mathcal{J}\right)$ into the uniform convergence
space $\left(Y,\mathcal{J}'\right)$.
\begin{eqnarray}
&\mbox{(1) }&\mbox{Then $f$ is uniformly continuous if
$\left(f\times
f\right)\left(\mathcal{U}\right)\in\mathcal{J}'$ whenever $\mathcal{U}\in\mathcal{J}$}\nonumber\\
&\mbox{(2) }&\mbox{The mapping $f$ is a uniform homeomorphism if
$f$ is bijective and
uniformly}\nonumber\\
&&\mbox{continuous with a uniformly continuous inverse.}\nonumber\\
&\mbox{(3) }&\mbox{The mapping $f$ is an embedding if it is a
uniform
homeomorphism onto its}\nonumber\\
&&\mbox{range, that is, $f: X\rightarrow f\left(X\right)\subseteq
Y$is a homeomorphism when $f\left(X\right)\subseteq Y$ is}\nonumber\\
&&\mbox{considered with the subspace uniform convergence
structure.}\nonumber
\end{eqnarray}
\end{definition}

With the concept of a uniform convergence space comes that of
Cauchy filter, completeness and completion.
\begin{definition}\label{CaychyComp}
Let $X$ be a uniform convergence space.
\begin{eqnarray}
&\mbox{(1) }&\mbox{A filter $\mathcal{F}$ on $X$ is a Cauchy
filter if $\mathcal{F}\times
\mathcal{F}\in\mathcal{J}$.}\nonumber\\
&\mbox{(2) }&\mbox{$X$ is complete if every Cauchy filter on $X$
converges.}\nonumber\\
&\mbox{(3) }&\mbox{A uniform convergence space $\widetilde{X}$ is
a completion of $X$ if there is a uniformly
continuous}\nonumber\\
&&\mbox{embedding $\iota:X\rightarrow\widetilde{X}$ such that
$\iota\left(X\right)$ is dense in $\widetilde{X}$ and, for every
complete uniform}\nonumber\\
&&\mbox{convergence space $Y$ and uniformly continuous mapping
$f:X\rightarrow Y$ there is a unique}\nonumber\\
&&\mbox{uniformly continuous mapping $f^{*}$ such that the diagram
commutes.}\nonumber
\end{eqnarray}\\
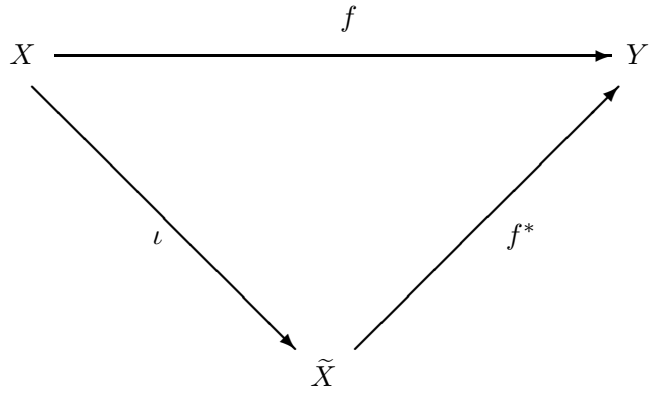
\begin{math}
\setlength{\unitlength}{1cm} \thicklines
\begin{picture}(13,6)

\put(3.4,5.4){$X$} \put(4.0,5.5){\vector(1,0){7.4}}
\put(11.6,5.4){$Y$} \put(7.8,5.9){$f$}
\put(3.7,5.1){\vector(1,-1){3.5}} \put(8.0,1.6){\vector(1,1){3.5}}
\put(7.4,1.1){$\widetilde{X}$} \put(5.3,3){$\iota$}
\put(10.0,3){$f^{*}$}

\end{picture}
\end{math}
\end{definition}
The main result concerning completions is the following
\cite{Gahler 2}.
\begin{theorem}\label{UCSCompletion}
Let $X$ be a Hausdorff uniform convergence space.  Then there
exists a uniform convergence space $\widetilde{X}$ that satisfies
Definition \ref{CaychyComp} (3).  Moreover, $\widetilde{X}$ is
unique up to homeomorphism.
\end{theorem}
The following result, found in \cite{Gahler 2} and \cite{Beattie},
concerns the completeness of a product space.
\begin{theorem}\label{ProdComplete}
Let $\left(X_{i}\right)_{i\in I}$ be a family of complete uniform
convergence spaces.  Then the product uniform convergence
structure on $\prod_{i\in I}X_{i}$ is complete.
\end{theorem}

We now prove three results concerning completions.
\begin{theorem}\label{ProdCom}
Let $\left(X_{i}\right)_{i\in I}$ be a family of uniform
convergence spaces and let $X$ denote their Cartesian product
equipped with the product uniform convergence structure.  Then the
completion $\widetilde{X}$ of $X$ is the product of the
completions $\widetilde{X}_{i}$ of the $X_{i}$.
\end{theorem}
\begin{proof}
First note that $\prod_{i\in I}X_{i}$ is complete by Theorem
\ref{ProdComplete}.  For every $i$, let
$\varphi_{i}:X_{i}\rightarrow \widetilde{X}_{i}$ be the uniformly
continuous embedding associated with the completion
$\widetilde{X}_{i}$ of $X_{i}$. Define the mapping
$\varphi:X\rightarrow \prod\widetilde{X}_{i}$ through
\begin{eqnarray}\label{VarphiDef}
\varphi:x=\left(x_{i}\right)\mapsto
\left(\varphi_{i}\left(x_{i}\right)\right)
\end{eqnarray}
Moreover, for each $i$, let $\pi_{i}:X\rightarrow X_{i}$ be the
projection mapping.\\
Since each $\varphi_{i}$ is injective, so is $\varphi$.  Moreover,
we have
\begin{eqnarray}
&\mathcal{U}\in\mathcal{J}_{P}&\Rightarrow\left(\pi_{i}\times
\pi_{i}\right)\left(\mathcal{U}\right)\in\mathcal{J}_{i}\nonumber\\
&& \Rightarrow\left(\varphi_{i}\times\varphi_{i}\right)
\left(\left(\pi_{i}\times\pi_{i}\right)\left(\mathcal{U}\right)\right)
\in\widetilde{\mathcal{J}}_{i}\nonumber\\
&&\Rightarrow\prod_{i\in
I}\left(\varphi_{i}\times\varphi_{i}\right)
\left(\left(\pi_{i}\times\pi_{i}\right)\left(\mathcal{U}\right)\right)
\in\widetilde{\mathcal{J}}_{P}\nonumber\\
&&\Rightarrow\left(\varphi\times\varphi\right)\left(\mathcal{U}\right)
\in\widetilde{\mathcal{J}}_{P}\nonumber
\end{eqnarray}
Hence $\varphi$ is uniformly continuous.  Similarly, if the filter
$\mathcal{V}$ on $\varphi\left(X\right)\times
\varphi\left(X\right)$ belongs to the subspace uniform convergence
structure, then
\begin{eqnarray}
&\left(\pi_{i}\times\pi_{i}\right)\left(\mathcal{V}\right)\in
\widetilde{\mathcal{J}}_{i}&\Rightarrow\left(\varphi_{i}^{-1}
\times\varphi_{i}^{-1}\right)\left(\left(\pi_{i}\times\pi_{i}\right)
\left(\mathcal{V}\right)\right)\in\mathcal{J}_{i}\nonumber\\
&&\Rightarrow \prod_{i\in I}\left(\varphi_{i}^{-1}
\times\varphi_{i}^{-1}\right)\left(\left(\pi_{i}\times\pi_{i}\right)
\left(\mathcal{V}\right)\right)\in\mathcal{J}_{P}\nonumber\\
&&\Rightarrow\left(\varphi^{-1}\times\varphi^{-1}\right)
\left(\mathcal{V}\right)\in\mathcal{J}_{P}\nonumber
\end{eqnarray}
so that $\varphi^{-1}$ is uniformly continuous.  Hence $\varphi$
is a uniformly continuous embedding.\\
The denseness of $\varphi\left(X\right)$ in $\prod_{i\in
I}\widetilde{X}_{i}$ follows by the denseness of
$\varphi_{i}\left(X_{i}\right)$ in $\widetilde{X}_{i}$, for each
$i\in I$.  The extension property of uniformly continuous mappings
into a complete uniform convergence space follows in the standard
way.
\end{proof}
\begin{theorem}\label{SubspUCSComp}
Let $X$ be a subspace of the uniform convergence space $Y$.  Then
there is a subspace $\widetilde{X}$ of the completion
$\widetilde{Y}$ of $Y$ such that $\widetilde{X}$ is the completion
of $X$.
\end{theorem}
\begin{proof}
Let $\varphi:Y\rightarrow \widetilde{Y}$ be the uniformly
continuous embedding associated with the completion
$\widetilde{Y}$ of $Y$.  Set
\begin{eqnarray}
\widetilde{X}=\{\widetilde{y}\in\widetilde{Y}\mbox{ :
}\exists\mbox{ }\mathcal{F}\mbox{ on }X\mbox{,
}\varphi\left(\mathcal{F}\right)\mbox{ converges to
}\widetilde{y}\}\label{CompDef}
\end{eqnarray}
Clearly $\varphi\left(X\right)$ is dense in $\widetilde{X}$.  To
see that $\widetilde{X}$ is complete, it suffices to show that it
is a closed subspace of $\widetilde{Y}$.  So let the filter
$\mathcal{G}$, with a trace on $\widetilde{X}$, converge to
$\widetilde{y}\in\widetilde{Y}$.  The the filter
\begin{eqnarray}
\mathcal{F}=[\{\varphi^{-1}\left(G\cap\widetilde{X}\right)\mbox{ :
}G\in\mathcal{G}\}]\nonumber
\end{eqnarray}
satisfies
\begin{eqnarray}
\varphi\left(\mathcal{G}\right)\mbox{ converges to }\widetilde{y}
\end{eqnarray}
so that $\widetilde{y}\in\widetilde{X}$.  Hence $\widetilde{X}$ is
a closed subspace of $\widetilde{Y}$.\\
The extension property for uniformly continuous mappings
$f:X\rightarrow Z$, where $Z$ is a complete uniform convergence
space, follows in the standard way.
\end{proof}\\
Lastly, we have the following.
\begin{theorem}\label{EmbExt}
Let $X$ and $Y$ be uniform convergence spaces, and $T:X\rightarrow
Y$ a uniformly continuous embedding.  Then there exists a
uniformly continuous embedding
$\widetilde{T}:\widetilde{X}\rightarrow\widetilde{Y}$, where
$\widetilde{X}$ and $\widetilde{Y}$ are the completions of $X$ and
$Y$ respectively, which extends $T$.
\end{theorem}
\begin{proof}
Since $X$ is homeomorphic to the subspace $T\left(X\right)$ of
$Y$, it follows that $\widetilde{X}$ is homeomorphic to
$\widetilde{T\left(X\right)}$ which is a subspace of
$\widetilde{Y}$ by Theorem \ref{SubspUCSComp}.  Hence we obtain
the following two commutative
diagrams.\\
\\
\begin{math}
\setlength{\unitlength}{1cm} \thicklines
\begin{picture}(13,6)

\put(3.4,5.4){$X$} \put(4.0,5.5){\vector(1,0){6.4}}
\put(10.6,5.4){$Y$} \put(6.8,5.7){$T$}
\put(3.5,5.2){\vector(0,-1){3.5}} \put(4.0,1.2){\vector(1,0){6.4}}
\put(3.4,1.1){$\widetilde{X}$} \put(10.6,1.1){$\widetilde{Y}$}
\put(3.2,3.4){$\phi$} \put(10.8,3.4){$\varphi$}
\put(10.7,5.2){\vector(0,-1){3.5}} \put(6.8,1.4){$\widetilde{T}$}

\end{picture}
\end{math}\\
\\
\begin{math}
\setlength{\unitlength}{1cm} \thicklines
\begin{picture}(13,6)

\put(3.4,5.4){$T\left(X\right)$} \put(4.4,5.5){\vector(1,0){6.0}}
\put(10.6,5.4){$Y$} \put(6.8,5.7){$i_{Y}$}
\put(3.5,5.2){\vector(0,-1){3.5}} \put(4.4,1.2){\vector(1,0){6.0}}
\put(3.4,1.1){$\widetilde{T\left(X\right)}$}
\put(10.6,1.1){$\widetilde{Y}$}
\put(2.4,3.4){$\varphi_{|T\left(X\right)}$}
\put(10.8,3.4){$\varphi$} \put(10.7,5.2){\vector(0,-1){3.5}}
\put(6.8,1.4){$i_{\widetilde{Y}}$}

\end{picture}
\end{math}\\
\\
Here $i_{Y}$ and $i_{\widetilde{Y}}$ denotes the inclusion
mappings on $Y$ and $\widetilde{Y}$, respectively.  It now follows
that $\widetilde{T}:\widetilde{X}\rightarrow \widetilde{Y}$ is an
embedding, which completes the proof.
\end{proof}

\end{article}
\end{document}